\RequirePackage{fix-cm}
\documentclass[smallextended]{svjour3}       
\smartqed  
\usepackage{graphicx}

\usepackage{algorithm}
\usepackage{algorithmic}

\begin{document}

\title{Adaptive iterative singular value thresholding algorithm to low-rank matrix recovery}


\author{Angang Cui$^{1,2}$\and
        Jigen Peng$^{3}$\and
        Haiyang Li$^{3}$.
}


\institute{
           1 School of Mathematics and Statistics, Xi'an Jiaotong University, Xi'an, 710049, China \\
           2 School of Mathematics and Statistics, Yulin University, Yulin, 719000, China\\
           3 School of Mathematics and Information Science, Guangzhou University, Guangzhou, 510006, China
}

\date{Received: date / Accepted: date}

\maketitle

\begin{abstract}
The problem of recovering a low-rank matrix from the linear constraints, known as affine matrix rank minimization problem, has been attracting extensive attention
in recent years. In general, affine matrix rank minimization problem is a NP-hard. In our latest work, a non-convex fraction
function is studied to approximate the rank function in affine matrix rank minimization problem and translate the NP-hard affine matrix rank minimization problem into a
transformed affine matrix rank minimization problem. A scheme of iterative singular value thresholding algorithm is generated to solve the regularized transformed affine matrix
rank minimization problem. However, one of the drawbacks for our iterative singular value thresholding algorithm is that the parameter $a$, which influences the behaviour of non-convex fraction function in the regularized transformed affine matrix rank minimization problem, needs to be determined manually in every simulation. In fact, how to determine the optimal parameter
$a$ is not an easy problem. Here instead, in this paper, we will generate an adaptive iterative singular value thresholding algorithm to solve the regularized transformed
affine matrix rank minimization problem. When doing so, our new algorithm will be intelligent both for the choice of the regularized parameter $\lambda$ and the parameter $a$. 

\end{abstract}

\section{Introduction}\label{intro}

The problem of recovering a low-rank or approximately low-rank matrix from the linear constraints, known as affine matrix rank minimization (AMRM) problem, has been attracting
extensive attention in recent years. Many applications such as minimum order system and low-dimensional Euclidean embedding in control theory \cite{Faz1,Faz2},
and collaborative filtering in recommender systems \cite{Cand3,Jan4} can be captured by solving the problem (AMRM). In mathematics, the problem (AMRM) can be described as the following
minimization:
\begin{equation}\label{equ1}
(\mathrm{AMRM})\ \ \ \ \ \ \min_{X\in R^{m\times n}} \ \mbox{rank}(X)\ \ s.t. \ \  \mathcal{A}(X)=b
\end{equation}
where $\mathcal{A}: R^{m\times n}\mapsto R^{d}$ is the linear map and the vector $b\in R^{d}$. Without loss of generality, in this paper, we assume $m\leq n$.
One important special case of the problem (AMRM) is the matrix completion (MC) problem \cite{Cand3,liu5,dong6,hu7,yu8,ma9,recht10}:
\begin{equation}\label{equ2}
(\mathrm{MC})\ \ \ \ \ \ \min_{X\in R^{m\times n}} \ \mbox{rank}(X)\ \ s.t. \ \  X_{i,j}=M_{i,j},\ \ (i,j)\in \Omega,
\end{equation}
which has been widely applied in famous Netflix problem, image inpainting problem, and so on. However, the problem (AMRM) is a challenging non-convex optimization
problem and is known as NP-hard \cite{recht10}.

Motivated by the recent development of non-convex relaxation approach in sparse signal recovery problems \cite{Xu11,Cao12,Zuo13}, in our latest work \cite{Cui14}, a continuous promoting low-rank
non-convex function
\begin{equation}\label{equ3}
P_{a}(X)=\sum_{i=1}^{m}\rho_{a}(\sigma_{i}(X))=\sum_{i=1}^{m}\frac{a\sigma_{i}(X)}{a\sigma_{i}(X)+1}
\end{equation}
in terms of the singular values of matrix $X$ is considered to substitute the rank function $\mathrm{rank}(X)$ in the NP-hard problem (AMRM), where $\sigma_{i}(X)$ represents the $i$-the largest 
singular value of matrix $X$, and the non-convex function
\begin{equation}\label{equ4}
\rho_{a}(t)=\frac{a|t|}{a|t|+1}\ \ \ (a>0)
\end{equation}
is the fraction function. It is clear to see that the non-convex function $P_{a}(X)$ has the rank approximation property \cite{Cui14}, with the change of parameter $a>0$,
it approximates the rank of matrix $X$:
\begin{equation}\label{equ5}
\displaystyle\lim_{a\rightarrow+\infty}P_{a}(X)=\displaystyle\lim_{a\rightarrow+\infty}\sum_{i=1}^{m}\frac{a\sigma_{i}(X)}{a\sigma_{i}(X)+1}
=\left\{
\begin{array}{ll}
0, & {\mathrm{if} \ \sigma_{i}(X)=0;} \\
\mathrm{rank}(X), & {\mathrm{if} \ \sigma_{i}(X)> 0.}
\end{array}
\right.
\end{equation}
Thus, by this transformation, we finally relax the NP-hard problem (AMRM) into the following transformed affine matrix rank minimization (TrAMRM) problem:
\begin{equation}\label{equ6}
(\mathrm{TrAMRM})\ \ \ \ \ \ \ \ \min_{X\in R^{m\times n}} \ P_{a}(X)\ \ s.t. \ \  \mathcal{A}(X)=b,
\end{equation}
where the non-convex surrogate function $P_{a}(X)$ is defined in (\ref{equ3}). 

Unfortunately, although we relax the NP-hard problem (AMRM) into a continuous problem
(TrAMRM), this relaxed problem is still computationally harder to solve due to the non-convex nature of the function $P_{a}(X)$, in fact, it is also NP-hrad. In \cite{Cui14},
we considered its regularized version:
\begin{equation}\label{equ7}
(\mathrm{RTrAMRM})\ \ \ \ \ \min_{X\in R^{m\times n}} \Big\{\|\mathcal{A}(X)-b\|_{2}^{2}+\lambda P_{a}(X)\Big\}
\end{equation}
where $\lambda>0$ is the regularized parameter.

As the unconstrained form, the problem (RTrAMRM) possesses much more algorithmic advantages. One nice property is that the proximal operator of fraction function has closed form analytical solutions for all values of parameter $a$. A scheme of iterative singular value thresholding algorithm (called ISVTA-Scheme 2 in \cite{Cui14}) has been devised to solve the problem (RTrAMRM) in our latest work \cite{Cui14}. A large number of numerical experiments have shown that the ISVTA-Scheme 2 can recover a low-rank matrix very well; however, we find that the parameter $a$, which influences the behaviour of non-convex fraction
function $\rho_{a}$ in ISVTA-Scheme 2, needs to be determined manually in every simulation. In fact, how to determine the optimal parameter $a$ in every
simulation is also a very hard problem. Unlike previous proposed ISVTA-Scheme 2 where the parameter $a$ needs to be determined manually in every
simulation, in this paper, we will generate an adaptive iterative singular value thresholding  algorithm (AISVTA) to solve the problem (RTrAMRM) which is intelligent both for the choice of the regularized
parameter $\lambda$ and the parameter $a$.

The rest of this paper is organized as follows. In Section \ref{section2}, we first review some known results from our latest work \cite{Cui14} for our previous proposed
ISVTA-Scheme 2, and then generate the AISVTA to solve the problem (RTrAMRM). In Section \ref{section3}, we test our algorithm on an image inpainting problem.
Finally, we give some concluding remarks in Section \ref{section4}.

\section{Algorithms for solving the problem (RTrAMRM)} \label{section2}

In this section, we first review some known results from \cite{Cui14} for our previous proposed ISVTA-Scheme 2, and then generate our AISVTA to solve the problem (RTrAMRM).
Unlike our previous proposed ISVTA-Scheme 2 where the parameter $a$ needs to be determined manually in every simulation, our newly proposed AISVTA will be intelligent both for the choice of the regularized parameter $\lambda$ and the parameter $a$.

\subsection{Iterative singular value thresholding algorithm (ISVTA)} \label{subsection2-1}

Define the proximal mapping of the non-convex function $P_{a}(X)$:
\begin{equation}\label{equ8}
 G_{a,\lambda}(Y):=\arg\min_{X\in R^{m\times n}}\Big\{\|X-Y\|_{F}^{2}+\lambda P_{a}(X)\Big\},
\end{equation}
we can get the following crucial result.

\begin{lemma}\label{lem1}
	Let $Y=U[\mathrm{Diag}(\sigma(Y),O_{m,n-m})]V^{\top}$ be the singular value decomposition (SVD) of matrix $Y\in R^{m\times n}$, where $O_{m,n-m}\in R^{m\times(n-m)}$ is a
	$m\times(n-m)$ zero matrix. Then the proximal operator $G_{a,\lambda}(Y)$ defined in (\ref{equ8})
	can be expressed as
	\begin{equation}\label{equ9}
    G_{a,\lambda}(Y)
    =U[\mathrm{Diag}(H_{a,\lambda}(\sigma(Y)),O_{m,n-m})]V^{\top}
	\end{equation}
	where 
	$$	H_{a,\lambda}(\sigma(Y)):=\big(h_{a,\lambda}(\sigma_{1}(Y),h_{a,\lambda}(\sigma_{2}(Y)),\cdots,h_{a,\lambda}(\sigma_{m}(Y))\big)^{\top}$$
	with 
	\begin{equation}\label{equ10}
    h_{a,\lambda}(\gamma)=\left\{
	\begin{array}{ll}
	g_{a,\lambda}(\gamma), & \ \ \mathrm{if} \ {|\gamma|> t_{a,\lambda};} \\\\
	0, & \ \ \mathrm{if} \ {|\gamma|\leq t_{a,\lambda}.}
	\end{array}
	\right.
	\end{equation}
	\begin{equation}\label{equ11}
	g_{a,\lambda}(\gamma)=\bigg(\frac{1+a|\gamma|}{3a}\Big(1+2\cos\big(\frac{\phi_{a,\lambda}(\gamma)}{3}-\frac{\pi}{3}\big)\Big)-\frac{1}{a}\bigg)\cdot sign(\gamma),
	\end{equation}
	$$\phi_{a,\lambda}(\gamma)=\arccos\Big(\frac{27\lambda a^{2}}{4(1+a|\gamma|)^{3}}-1\Big),$$
	\begin{equation}\label{equ12}
	t_{a,\lambda}=\left\{
	\begin{array}{ll}
	\frac{\lambda a}{2}, & \ \ \mathrm{if} \ {\lambda\leq \frac{1}{a^{2}};} \\\\
	\sqrt{\lambda}-\frac{1}{2a}, & \ \ \mathrm{if} \ {\lambda>\frac{1}{a^{2}}.}
	\end{array}
	\right.
	\end{equation}
\end{lemma}

In the following, the ISVTA is generated to solve the problem (RTrAMRM). We consider the following regularization function
\begin{equation}\label{equ13}
\mathcal{C}_{\lambda}(X)=\|\mathcal{A}(X)-b\|_{2}^{2}+\lambda P_{a}(X)
\end{equation}
and its surrogate function
\begin{equation}\label{equ14}
\mathcal{C}_{\lambda, \mu}(X, Z)=\mu[\mathcal{C}_{\lambda}(X)-\|\mathcal{A}(X)-\mathcal{A}(Z)\|_{2}^{2}]+\|X-Z\|_{F}^{2}
\end{equation}
for any fixed $\lambda>0$, $\mu>0$ and $Z\in R^{m\times n}$. When we set $0<\mu\leq \frac{1}{\|\mathcal{A}\|_{2}^{2}}$, we can get that
$$\|X-Z\|_{F}^{2}-\mu\|\mathcal{A}(X)-\mathcal{A}(Z)\|_{F}^{2}\geq0.$$
Therefore, we have
\begin{equation}\label{equ15}
\begin{array}{llll}
\mathcal{C}_{\lambda, \mu}(X, Z)&=&\mu\mathcal{C}_{\lambda}(X)+\|X-Z\|_{F}^{2}-\mu\|\mathcal{A}(X)-\mathcal{A}(Z)\|_{2}^{2}\\
&\geq&\mu\mathcal{C}_{\lambda}(X).
\end{array}
\end{equation}
Under the condition $0<\mu\leq\frac{1}{\|\mathcal{A}\|_{2}^{2}}$, if we suppose that the matrix $X^{\star}\in R^{m\times n}$ is a minimizer of the function
$\mathcal{C}_{\lambda}(X)$, then
\begin{eqnarray*}
	\mathcal{C}_{\lambda, \mu}(X,X^{\star})&=&\mu[\mathcal{C}_{\lambda}(X)-\|\mathcal{A}(X)-\mathcal{A}(X^{\star})\|_{2}^{2}]+\|X-X^{\star}\|_{F}^{2}\\
	&=&\mu \mathcal{C}_{\lambda}(X)\\
	&\geq&\mu \mathcal{C}_{\lambda}(X^{\star})\\
	&=&\mathcal{C}_{\mu}(X^{\star},X^{\star}),
\end{eqnarray*}
which implies that $X^{\star}$ is also a minimizer of $\mathcal{C}_{\lambda, \mu}(X,X^{\star})$. On the other hand, $\mathcal{C}_{\lambda, \mu}(X,Z)$ with $Z=X^{\star}$ can be reexpressed as
\begin{eqnarray*}
	\mathcal{C}_{\lambda,\mu}(X,X^{\star})
	&=&\|X-(X^{\star}-\mu \mathcal{A}^{\ast}\mathcal{A}(X^{\star})+\mu \mathcal{A}^{\ast}(b))\|_{F}^{2}\\
	&&+\lambda\mu P_{a}(X)+\mu\|b\|_{2}^{2}+\|X^{\star}\|_{F}^{2}-\mu\|\mathcal{A}(X^{\star})\|_{2}^{2}\\
	&&-\|X^{\star}-\mu \mathcal{A}^{\ast}\mathcal{A}(X^{\star})+\mu \mathcal{A}^{\ast}(b)\|_{F}^{2}\\
	&=&\|X-B_{\mu}(X^{\star})\|_{F}^{2}+\lambda\mu P_{a}(X)+\mu\|b\|_{2}^{2}\\
	&&+\|X^{\star}\|_{F}^{2}-\mu\|\mathcal{A}(X^{\star})\|_{2}^{2}-\|B_{\mu}(X^{\star})\|_{F}^{2},
\end{eqnarray*}
where $B_{\mu}(X^{\star})=X^{\star}+\mu \mathcal{A}^{\ast}(b-\mathcal{A}(X^{\star}))$. This implies that for any fixed $\lambda>0$ and $\mu>0$, minimizing
$\mathcal{C}_{\lambda,\mu}(X,X^{\star})$ on $X$ is equivalent to solve the following minimization problem 
\begin{equation}\label{equ16}
\min_{X\in R^{m\times n}}\Big\{\|X-B_{\mu}(X^{\star})\|_{F}^{2}+\lambda\mu P_{a}(X)\Big\}.
\end{equation}
By Lemma \ref{lem1}, the minimizer $X^{\star}$ of minimization problem  (\ref{equ16}) is given by
\begin{equation}\label{equ17}
\begin{array}{llll}
X^{\star}&=&G_{a,\lambda\mu}(B_{\mu}(X^{\star}))\\
&=&
U^{\star}[\mathrm{Diag}(H_{a,\lambda\mu}(\sigma(B_{\mu}(X^{\star}))),O_{m,n-m})](V^{\star})^{\top},
\end{array}
\end{equation}
where $U^{\star}[\mathrm{Diag}(\sigma(B_{\mu}(X^{\star})),O_{m,n-m})](V^{\star})^{\top}$ is the SVD of matrix $B_{\mu}(X^{\star})$, and  $G_{a,\lambda\mu}(\cdot)$ is obtained by replacing $\lambda$ with $\lambda\mu$ in $G_{a,\lambda}(\cdot)$. 

With the representation (\ref{equ17}), the ISVTA for solving the problem (RTrAMRM) can be naturally given by
\begin{equation}\label{equ18}
\begin{array}{llll}
X^{k+1}&=&G_{a,\lambda\mu}(B_{\mu}(X^{k}))\\
&=&
U^{k}[\mathrm{Diag}(H_{a,\lambda\mu}(\sigma(B_{\mu}(X^{k}))),O_{m,n-m})](V^{k})^{\top},
\end{array}
\end{equation}
where $U^{k}[\mathrm{Diag}(\sigma(B_{\mu}(X^{k})),O_{m,n-m})](V^{k})^{\top}$ is the SVD of matrix $B_{\mu}(X^{k})$. 

The basic convergence theorem of iteration (\ref{equ18}) can be stated as below.
\begin{theorem} \label{the1} {\rm(see \cite{Cui14})}
	Let $\{X^{k}\}$ be the sequence generated by the iteration (\ref{equ18}) with the step size $\mu$ satisfying $0<\mu<\frac{1}{\|\mathcal{A}\|_{2}^{2}}$. Then
	\begin{description}
		\item[1)] The sequence $\mathcal{C}_{\lambda}(X^{k})$ is decreasing;
		\item[2)] $\{X^{k}\}$ is asymptotically regular, i.e., $\lim_{k\rightarrow\infty}\|X^{k+1}-X^{k}\|_{F}^{2}=0$;
		\item[3)] Any accumulation point of $\{X^{k}\}$ is a stationary point of the problem (RTrAMRM).
	\end{description}
\end{theorem}

According to iteration $(\ref{equ18})$, in \cite{Cui14}, two schemes of ISVTA (ISVTA-Scheme 1 and ISVTA-Scheme 2) are generated to solve the problem (RTrAMRM). Especially, in ISVTA-Scheme 2\cite{Cui14,Cui14jia}, and adaptive strategy is accepted to select the proper regularized parameter $\lambda$. 

Suppose that the matrix $X^{\star}$ of rank $r$ is
the optimal solution to the problem (RTrAMRM). In each iteration, the regularized parameter $\lambda$ can be selected as
\begin{equation}\label{equ19}
\begin{array}{llll}
\lambda=\left\{
\begin{array}{ll}
\lambda_{1,k}=\frac{2\sigma_{r+1}(B_{\mu}(X^{k}))}{a\mu}, & {\mathrm{if}\  \lambda_{1,k}\leq\frac{1}{a^{2}\mu};} \\\\
\lambda_{2,k}=\frac{(1-\xi)(2a\sigma_{r}(B_{\mu}(X^{k}))+1)^{2}}{4a^{2}\mu}, & {\mathrm{if}\  \lambda_{1,k}>\frac{1}{a^{2}\mu}.}
\end{array}
\right.
\end{array}
\end{equation}
where $\xi>0$\ is a very small small positive number such as 0.01 or 0.001. Using the adaptive strategy (\ref{equ19}), the ISVTA-Scheme 2 will be adaptive for the choice of the regularization parameter $\lambda$ in each iteration. The ISVTA-Scheme 2 is summarized in Algorithm \ref{alg:A}.

\begin{algorithm}[h!]
	\caption{: ISVTA-Scheme 2}
	\label{alg:A}
	\begin{algorithmic}
		\STATE {\textbf{Input}: $\mathcal{A}: R^{m\times n}\mapsto R^{d}$, $b\in R^{d}$,  $\mu\in(0,\frac{1}{\|\mathcal{A}\|_{2}^{2}})$, $a=a_{0}$ ($a_{0}$ is a given positive number), $\xi>0$ is a very small positive number such as 0.01 or 0.001, $k=0$;}
		\STATE {\textbf{Initialize}: Given $X^{0}\in R^{m\times n}$;}
		\STATE {\textbf{while} not converged \textbf{do}}
		\STATE \ \ \ \ {$B_{\mu}(X^{k})=X^{k}+\mu \mathcal{A}^{\ast}(b-\mathcal{A}(X^{k}))$;}
		\STATE \ \ \ \ {Compute the SVD of $B_{\mu}(X^{k})$ as: $B_{\mu}(X^{k})=U^{k}[\mathrm{Diag}(\sigma(B_{\mu}(X^{k}))),O_{m\times(n-m)}](V^{k})^{\top}$;}
		\STATE \ \ \ \ {$\lambda_{1,k}=\frac{2\sigma_{r+1}(B_{\mu}(X^{k}))}{a\mu}$, $\lambda_{2,k}=\frac{(1-\xi)(2a\sigma_{r}(B_{\mu}(X^{k}))+1)^{2}}{4a^{2}\mu}$;}
		\STATE \ \ \ \ {if\ $\lambda_{1,k}\leq\frac{1}{a^{2}\mu}$\ then}
		\STATE \ \ \ \ \ \ \ {$\lambda=\lambda_{1,k}$, $t_{a,\lambda\mu}=\frac{\lambda\mu a}{2}$}
		\STATE \ \ \ \ \ \ \ {for\ $i=1:m$}
		\STATE \ \ \ \ \ \ \ \ \ {1.\ $\sigma_{i}(B_{\mu}(X^{k}))>t_{a,\lambda\mu}$, then $\bar{\sigma}_{i}=g_{a,\lambda\mu}(\sigma_{i}(B_{\mu}(X^{k})))$;}
		\STATE \ \ \ \ \ \ \ \ \ {2.\ $\sigma_{i}(B_{\mu}(X^{k}))\leq t_{a,\lambda\mu}$, then $\bar{\sigma}_{i}=0$;}
		\STATE \ \ \ \ \ \ \ {end}
		\STATE \ \ \ \ {else}
		\STATE \ \ \ \ \ \ \ {$\lambda=\lambda_{2,k}$, $t_{a,\lambda\mu}=\sqrt{\lambda\mu}-\frac{1}{2a}$;}
		\STATE \ \ \ \ \ \ \ {for\ $i=1:m$}
		\STATE \ \ \ \ \ \ \ \ \ \ {1.\ $\sigma_{i}(B_{\mu}(X^{k}))>t_{a,\lambda\mu}$, then $\bar{\sigma}_{i}=g_{a,\lambda\mu}(\sigma_{i}(B_{\mu}(X^{k})))$;}
		\STATE \ \ \ \ \ \ \ \ \ \ {2.\ $\sigma_{i}(B_{\mu}(X^{k}))\leq t_{a,\lambda\mu}$, then $\bar{\sigma}_{i}=0$;}
		\STATE \ \ \ \ \ \ \ {end}
		\STATE \ \ \ \ {end}
		\STATE \ \ \ \ \ \ \ \ \ {$X^{k+1}=U^{k}[\mathrm{Diag}(\bar{\sigma}),O_{m\times(n-m)}](V^{k})^{\top}$;}
		\STATE \ \ \ \ \ \ \ \ \ {$k\rightarrow k+1$;}
		\STATE{\textbf{end while}}
		\STATE{\textbf{return}: $X^{opt}$}
	\end{algorithmic}
\end{algorithm}

A large number of numerical experiments on some completion of low-rank random matrices and image inpainting problems have shown that the ISVTA-Scheme 2 performances very
well in recovering a low-rank matrix compared with some state-of-art methods. One of the drawbacks for our ISVTA-Scheme 2 is that the parameter $a$, which influences the behaviour of
non-convex fraction function $\rho_{a}$, needs to be determined manually in every simulation. In fact, how to determine the best parameter $a$ is not an easy problem.

\subsection{Adaptive iterative singular value thresholding algorithm (AISVTA)} \label{subsection2-2}

Different from our previous proposed ISVTA-Scheme 2 where the parameter $a$ needs to be determined manually in every simulation, in this section, we will generate an adaptive iterative singular value thresholding algorithm (AISVTA) to solve the problem (RTrAMRM). AISVTA will be intelligent both for the choice of the regularized parameter $\lambda$ and the parameter $a$, which is one of the advantages for the AISVTA compared with the our previous proposed ISVTA-Scheme 2. 

In the following descriptions, we will generate our AISVTA to solve the problem (RTrAMRM).

\begin{lemma}\label{lem2}
\cite{Cui15} For any $\beta\geq0$, $\gamma\in R$ and $0<a\leq\frac{1}{\sqrt{\lambda}}$, the function 
\begin{equation}\label{equ20}
f_{a,\lambda}(\beta):=(\beta-\gamma)^{2}+\lambda\frac{a\beta}{a\beta+1}
\end{equation}
is strictly convex. 
\end{lemma}

Lemma \ref{lem2} told us that, when the parameter satisfies  $0<a\leq\frac{1}{\sqrt{\lambda}}$, the function 
$f_{a,\lambda}(\beta)$ defined in (\ref{equ20}) is a strictly convex function, which implies that there exist the unique minimizer for the function 
$f_{a,\lambda}(\beta)$ defined in (\ref{equ20}). The following theorem gives the expression of this unique minimizer of the function $f_{a,\lambda}(\beta)$  defined in (\ref{equ20}).
\begin{theorem}\label{theorem1}
	\cite{Cui15}For any fixed $\lambda>0$ and  $0<a\leq\frac{1}{\sqrt{\lambda}}$, suppose $\beta_{\lambda}$ is the minimizer of the minimization problem 
	\begin{equation}\label{equ21}
	\min_{\beta\geq0}\Big\{(\beta-\gamma)^{2}+\lambda\frac{a\beta}{a\beta+1}\Big\}, 
	\end{equation}
	then $\beta_{\lambda}$ is unique and 
	\begin{equation}\label{equ22}
	\beta_{\lambda}=\bar{h}_{a,\lambda}(\gamma)
	:=\left\{
	\begin{array}{ll}
	g_{a,\lambda}(\gamma), & \ \ \ {\gamma>\bar{t}_{a,\lambda};} \\
	0, & \ \ \ {\gamma\leq \bar{t}_{a,\lambda}.}
	\end{array}
	\right.
	\end{equation}
	where
	\begin{equation}\label{equ23}
	\bar{t}_{a,\lambda}=\frac{\lambda a}{2}
	\end{equation}
\end{theorem}
and $g_{a,\lambda}(\cdot)$ is defined in Lemma \ref{lem1}. 

Next, under the condition $0<a\leq\frac{1}{\sqrt{\lambda}}$, we shall present the AISVTA to solve the problem (RTrAMRM. 

Similar as the generation of iteration (\ref{equ18}), under the condition $0<a\leq\frac{1}{\sqrt{\lambda\mu}}$, the ISVTA for solving the problem (RTrAMRM) can be rewritten as
\begin{equation}\label{equ24}
X^{k+1}=
U^{k}[\mathrm{Diag}(\bar{H}_{a,\lambda\mu}(\sigma(B_{\mu}(X^{k}))),O_{m,n-m})](V^{k})^{\top},
\end{equation}
where 
\begin{eqnarray*}	&&\bar{H}_{a,\lambda\mu}(\sigma(B_{\mu}(X^{k})))\\
&&=\big(\bar{h}_{a,\lambda\mu}(\sigma_{1}(B_{\mu}(X^{k})),\bar{h}_{a,\lambda}(\sigma_{2}(B_{\mu}(X^{k}))),\cdots,\bar{h}_{a,\lambda}(\sigma_{m}(B_{\mu}(X^{k})))\big)^{\top},
\end{eqnarray*}
and $U^{k}[\mathrm{Diag}(\sigma(B_{\mu}(X^{k})),O_{m,n-m})](V^{k})^{\top}$ is the SVD of matrix $B_{\mu}(X^{k})$.

Here, we will
generate an adaptive rule for the choice of the parameters $\lambda$ and $a$ in our iteration (\ref{equ24}). When doing so, the iteration (\ref{equ24}) will
be intelligent both for the choice of the regularized parameter $\lambda$ and parameter $a$.

1) Adaptive for the choice of parameter $a$: Note that the parameter $a$ in iteration (\ref{equ24}) should be satisfied $0<a\leq\frac{1}{\sqrt{\lambda\mu}}$. Therefore, we can choose the parameter a as
\begin{equation}\label{equ25}
a=\frac{\tau}{\sqrt{\lambda\mu}},
\end{equation}
where $\tau\in(0,1]$ is a given positive number. When we 
set $a=\frac{\tau}{\sqrt{\lambda\mu}}$, the threshold function $\bar{t}_{a,\lambda\mu}=\frac{\lambda\mu a}{2}$ can be rewritten as 
\begin{equation}\label{equ26}
\bar{t}_{a,\lambda\mu}=\frac{\tau\sqrt{\lambda\mu}}{2}.
\end{equation}
To see clear that once the value of the regularized  parameter $\lambda$ is determined, the parameter $a$ can be given by (\ref{equ25}), and therefore the iteration (\ref{equ24}) will be adaptive for the choice of the parameter $a$. For the choice of the proper regularized parameter $\lambda$, here, the rule which
is used to select the proper regularized parameter $\lambda$ in our previous proposed ISVTA-Scheme 2 is again used to select the proper regularized parameter $\lambda$ in iteration (\ref{equ24}).

2) Adaptive for the choice of regularized parameter $\lambda$: Let the matrix $X^{\star}$ of rank $r$ be the optimal solution to the  problem (RTrAMRM) under the condition $a=\frac{\tau}{\sqrt{\lambda\mu}}$. Then, the following inequalities hold
$$\sigma_{i}(B_{\mu}(X^{\star}))>\frac{\tau\sqrt{\lambda\mu}}{2}\Leftrightarrow i\in\{1,2,\cdots,r\},$$
$$\sigma_{j}(B_{\mu}(X^{\star}))\leq \frac{\tau\sqrt{\lambda\mu}}{2}\Leftrightarrow j\in\{r+1,r+2,\cdots,m\},$$
which implies that 
$$\frac{4\sigma^{2}_{r+1}(B_{\mu}(X^{\star}))}{\tau^{2}\mu}\leq\lambda<\frac{4\sigma^{2}_{r}(B_{\mu}(X^{\star}))}{\tau^{2}\mu}.$$
The above estimation provides an exact location of the regularized parameter $\lambda$. Here, we can choose the regularized parameter $\lambda$ as 
$$\lambda^{\star}=\frac{4}{\tau^{2}\mu}\Big[(1-\alpha)\sigma^{2}_{r+1}(B_{\mu}(X^{\star}))+\alpha\sigma^{2}_{r}(B_{\mu}(X^{\star})) \Big],$$
where $\alpha\in[0,1)$. When we set $\alpha=0$, a most reliable choice of the proper regularized parameter $\lambda$ specified by
\begin{equation}\label{equ27}
\lambda^{\star}=\frac{4}{\tau^{2}\mu}\sigma^{2}_{r+1}(B_{\mu}(X^{\star})).
\end{equation}
Combing with (\ref{equ25}) and (\ref{equ27}), we can get an adaptive strategy for the choice of regularized parameter $\lambda$ and parameter $a$ in iteration (\ref{equ24}): 
\begin{equation}\label{equ28}
\left\{
\begin{array}{ll}
\lambda^{\star}=\displaystyle\frac{4}{\tau^{2}\mu}\sigma^{2}_{r+1}(B_{\mu}(X^{\star})); \\\\
a^{\star}=\displaystyle\frac{\tau}{\sqrt{\lambda^{\star}\mu}}.
\end{array}
\right.
\end{equation}
In each iteration, we can approximate the optimal solution $X^{\star}$ by $X^{k}$. Then, in each iteration, the proper regularized parameter $\lambda$ and parameter $a$ in iteration (\ref{equ24}) can be selected as
\begin{equation}\label{equ29}
\left\{
\begin{array}{ll}
\lambda_{k}^{\star}=\displaystyle\frac{4}{\tau^{2}\mu}\sigma^{2}_{r+1}(B_{\mu}(X^{k})); \\\\
a_{k}^{\star}=\displaystyle\frac{\tau}{\sqrt{\lambda_{k}^{\star}\mu}},
\end{array}
\right.
\end{equation}
where $\tau\in(0,1]$.

\begin{algorithm}[h!]
	\caption{: AISVTA}
	\label{alg:B}
	\begin{algorithmic}
		\STATE {\textbf{Input}: $\mathcal{A}: R^{m\times n}\mapsto R^{d}$, $b\in R^{d}$,  $\mu\in(0,\frac{1}{\|\mathcal{A}\|_{2}^{2}})$, $\tau\in(0,1]$, $k=0$;}
		\STATE {\textbf{Initialize}: Given $X^{0}\in R^{m\times n}$;}
		\STATE {\textbf{while} not converged \textbf{do}}
		\STATE \ \ \ \ {$B_{\mu}(X^{k})=X^{k}+\mu \mathcal{A}^{\ast}(b-\mathcal{A}(X^{k}))$;}
		\STATE \ \ \ \ {Compute the SVD of $B_{\mu}(X^{k})$ as: $B_{\mu}(X^{k})=U^{k}[\mathrm{Diag}(\sigma(B_{\mu}(X^{k}))),O_{m\times(n-m)}](V^{k})^{\top}$;}
		\STATE \ \ \ \ {$\lambda_{k}^{\star}=\displaystyle\frac{4}{\tau^{2}\mu}\sigma^{2}_{r+1}(B_{\mu}(X^{k}))$;}
		\STATE \ \ \ \ {if\ $\lambda_{k}^{\star}\neq0$\ then}
		\STATE \ \ \ \ \ \ \ {$\lambda=\lambda_{k}^{\star}$, $a=a_{k}^{\star}=\displaystyle\frac{\tau}{\sqrt{\lambda_{k}^{\star}\mu}}$, $\bar{t}_{a,\lambda\mu}=\frac{\lambda\mu a}{2}$;}
		\STATE \ \ \ \ \ \ \ {for\ $i=1:m$}
		\STATE \ \ \ \ \ \ \ \ \ {1.\ $\sigma_{i}(B_{\mu}(X^{k}))>\bar{t}_{a,\lambda\mu}$, then $\bar{\sigma}_{i}=g_{a,\lambda\mu}(\sigma_{i}(B_{\mu}(X^{k})))$;}
		\STATE \ \ \ \ \ \ \ \ \ {2.\ $\sigma_{i}(B_{\mu}(X^{k}))\leq \bar{t}_{a,\lambda\mu}$, then $\bar{\sigma}_{i}=0$;}
		\STATE \ \ \ \ \ \ \ {end}
		\STATE \ \ \ \ \ \ \  {$X^{k+1}=U^{k}[\mathrm{Diag}(\bar{\sigma}),O_{m\times(n-m)}](V^{k})^{\top}$;}
		\STATE \ \ \ \ {else}
		\STATE \ \ \ \ \ \ \ {$X^{k+1}=B_{\mu}(X^{k}))$;}
		\STATE \ \ \ \ {end}
		\STATE \ \ \ \ {$k\rightarrow k+1$;}
		\STATE{\textbf{end while}}
		\STATE{\textbf{return}: $X^{opt}$}
	\end{algorithmic}
\end{algorithm}

It is important to note that, in some iterations, the value of $\lambda^{\star}_{k}$ may be 0. If $\lambda^{\star}_{k}=0$, by (\ref{equ29}), the value 
of $a_{k}^{\star}$ may be $+\infty$. In computer operations, we must try to avoid this situation. In fact, if $\lambda=\lambda^{\star}_{k}=0$ in some iterations, the minimization problem 
\begin{equation}\label{equ30}
\min_{X\in R^{m\times n}}\Big\{\|X-B_{\mu}(X^{k})\|_{F}^{2}+\lambda\mu P_{a}(X)\Big\}
\end{equation}
will reduces to 
\begin{equation}\label{equ31}
\min_{X\in R^{m\times n}}\|X-B_{\mu}(X^{k})\|_{F}^{2}. 
\end{equation}
Under this situation, the minimizer of (\ref{equ30}) can be written as 
\begin{equation}\label{equ32}
X^{k+1}=B_{\mu}(X^{k}).
\end{equation}
This implies that, during the iteration process, the situation $a_{k}^{\star}=+\infty$ can be completely avoided.

By above operations, the iteration (\ref{equ24}) will be adaptive for the choice of the regularized 
parameter $\lambda$ and parameter $a$ in each 
iteration. The iteration (\ref{equ24}) with the parameter choice strategy (\ref{equ29}) is our  AISVTA, and it is summarized in Algorithm \ref{alg:B}.

\section{Numerical experiments} \label{section3}

\begin{figure}[h!]
	\centering
	\begin{minipage}[t]{0.47\linewidth}
		\centering
		\includegraphics[width=1.1\textwidth]{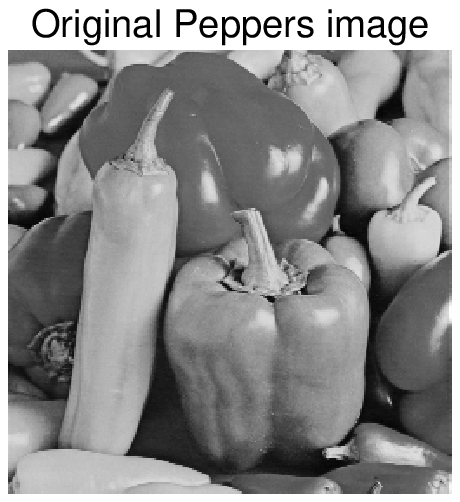}
	\end{minipage}
	\begin{minipage}[t]{0.47\linewidth}
		\centering
		\includegraphics[width=1.1\textwidth]{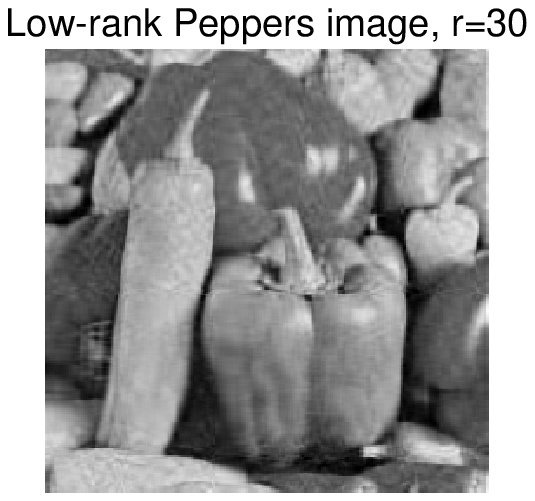}
	\end{minipage}
	\caption{Original $256\times 256$ Peppers image and its approximated low-rank image with rank $r=30$.} \label{figure1}
\end{figure}

In this section, to study the performance of the proposed AISVTA, some simulation experiments (for image impainting problem)
are considered. We compare our AISVTA with our previous proposed ISVTA-Scheme 2 on an image inpainting problem (Low-rank Peppers image inpainting problem) under the noise case. We test these two algorithms on a gray-scale
images: $256\times 256$ Peppers image. We use the SVD to obtain its approximated low-rank image with rank $r=30$. The original images, and its low-rank images are displayed in Figure \ref{figure1}.

We let $fr=s/r(m+n-r)$
to denote the freedom ration\cite{don16}, i.e., the ratio between the number of sampled entries and the `true dimensionality' of a $m\times n$ image of rank $r$, where $s$ represents the number of randomly sampled entries. If sampling ration is given as $sr\in[0,1]$, in these numerical experiments, the the number of randomly sampled entries $s$ can be obtained by using Matlab code: $s=\texttt {round}(sr*m*n)$. The stopping criterion is defined as
$$\frac{\|X^{k+1}-X^{k}\|_{F}}{\max\{1,\|X^{k}\|_{F}\}}\leq 10^{-8}$$
or maximal steps 5000, where $X^{k+1}$ and $X^{k}$ are numerical results from two continuous iterative steps.

For the given truth low-rank  $M\in R^{m\times n}$, we measure the accuracy of the generated solution $X^{opt}$ of our algorithms by the relative error ($\mathrm{RE}$)
$$\mathrm{RE}=\frac{\|X^{opt}-M\|_{F}}{\|M\|_{F}}.$$
In low-rank Peppers image impainting problem, the observed entries are polluted by noise: 
$$Q_{i,j}=M_{i,j}+\xi_{1}\ast E_{i,j}, \ (i,j)\in\Omega,$$
where $\xi_{1}\in(0,1)$, which means that  $\mathcal{P}_{\Omega}(Q)=\mathcal{P}_{\Omega}(M)+\mathcal{P}_{\Omega}(\xi_{1}\ast E)$. In 
these numerical experiments, the noise $E$ is generated by Matlab code: 
$$E=\texttt{randn}(m,n).$$
In AISVTA, we set $\tau=0.45$ and $\mu=0.99$. In ISVTA, we set $a=1$, $\mu=0.99$ and $\xi=0.01$. All the simulation experiments are performed using the Matlab R2015b on ThikPad S2 (Intel(R) Core(TM) i7-8565U CPU @ 1.80GHZ 1.99GHZ with 8GB of RAM running Microsoft Windows 10).

\begin{table}[h!]
\centering
\begin{tabular}{|c||l|l|l|l|}\hline
\multicolumn{5}{|c|}{SR=0.50}\\\hline
Peppers&\multicolumn{2}{|c}{AISVTA}&\multicolumn{2}{|c|}{ISVTA-Scheme 2}\\
\hline
($\xi_{1}$, $fr$)&RE&Time&RE&Time\\
\hline
(0.01, 2.2661)&1.56e-02&3.84&1.54e-02&4.73\\
\hline
(0.03, 2.2661)&4.88e-02&2.51&4.74e-02&7.40\\
\hline
(0.06, 2.2661)&9.21e-02&1.54&9.56e-02&11.20\\
\hline
\multicolumn{5}{|c|}{SR=0.40}\\\hline
Peppers&\multicolumn{2}{|c}{AISVTA}&\multicolumn{2}{|c|}{ISVTA-Scheme 2}\\
\hline
($\xi_{1}$, $fr$)&RE&Time&RE&Time\\
\hline
(0.01, 1.8129)&2.06e-02&11.15&2.05e-02&12.76\\
\hline
(0.03, 1.8129)&6.10e-02&6.01&6.67e-02&25.96\\
\hline
(0.06, 1.8129)&1.05e-01&3.42&1.43e-01& 71.04\\
\hline
\end{tabular}
\caption{\scriptsize Numerical results of ISVTA-Scheme 2 and AISVTA for low-rank Peppers image inpainting problems, $r=30$.}\label{table1}
\end{table}

Table \ref{table1} reports the numerical results of AISVTA and ISVTA-Scheme 2 for low-rank Peppers image inpainting problem with different $\mathrm{SR}$ and $\xi_{1}$. Comparing with these numerical results, we can see that the AISVTA and ISVTA-Scheme 2 
have almost the same recovery results, but the AISVTA has a faster running speed with the increasing of the value of $\xi_{1}$.

\section{Conclusions}\label{section4}
In this paper, we first review some known results from our lately work for ISVTA-Scheme 2, and then generate an AISVTA to solve the problem (RTrAMRM). Different from our previous proposed ISVTA-Scheme 2 where the parameter $a$ needs to be determined manually in every simulation, the AISVTA will be intelligent both for the choice of the regularized parameter $\lambda$ and the parameter $a$, which is one of the advantages for the AISVTA compared with our previous proposed ISVTA-Scheme 2. Numerical experiments on an image inpainting problem have shown that, under the noise case, the AISVTA and ISVTA-Scheme 2 have almost the same recovery results, but the AISVTA has a faster running speed with the increasing of the value of $\xi_{1}$.

\begin{acknowledgements}
The work was supported by
the National Natural Science Foundations of China (11771347, 91730306, 41390454, 11271297).
\end{acknowledgements}


\begin{thebibliography}{}

\bibitem{Faz1}
Fazel, M., Hindi, H., Boyd, S.P.: A rank minimization heuristic with application to minimum order system approximation,
\textit{In proceedings of American Control Conference, Arlington, VA}, 2001, 6, pp.~4734--4739

\bibitem{Faz2}
Fazel, M., Hindi, H., Boyd, S.P.: Log-det heuristic for matrix rank minimization with applications to Hankel and Euclidean distance matrices,
\textit{In Proceedings of American Control Conference, Denever, Colorado}, 2003, 3, pp.~2156--2162


\bibitem{Cand3}
Candes, E.J., Rechtc, B.: Exact matrix completion via convex optimization,\textit{Foundations of Computational Mathematics}, 2009, 9, pp.~717--772

\bibitem{Jan4}
Jannach, D., Zanker, M., Felfernig, A., Friedrich, G.: Recommender Systems: An Introduction, \textit{Cambridge University Press}, 2010


\bibitem{liu5}
Liu, P., Lewis, J., Rhee, T.: Low-rank matrix completion to reconstruct incomplete rendering images, \textit{IEEE Transactions on Visualization and Computer Graphics Volume}, 2018,
24(8), pp.~2353--2365


\bibitem{dong6}
Dong, J., Xue, Z. C., Guan, J., Han, Z. F., Wang, W. W.: Low rank matrix completion using truncated nuclear norm and sparse regularizer, \textit{Signal Processing: Image Communication},
2018, 68, pp.~76--87



\bibitem{hu7}
Hu, Y., Zhang, D., Ye, J., Li, X., He, X.: Fast and accurate matrix completion via truncated nuclear norm regularization, \textit{IEEE Transactions on Pattern Analysis and Machine Intelligence}, 2013, 35(9), pp.~2117--2130


\bibitem{yu8}
Yu, Y.C., Peng,J.G.: A modified primal-dual method with applications to some sparse recovery problems, \textit{Applied mathematics and computation}, 2018, 333, pp.~76--94


\bibitem{ma9}
Ma, S., Goldfarb, D., Chen, L.: Fixed point and Bregman iterative methods for matrix rank minimization, \textit{Mathematical Programming, Ser. A}, 2011, 128, pp.~321--353

\bibitem{recht10}
Recht, B., Fazel, M., Parrilo, P.A.: Guaranteed minimum-rank solution of linear matrix equations via nuclear norm minimization, \textit{SIAM Review}, 2010, 52, pp.~471--501

\bibitem{Xu11}
Xu, Z.B., Chang, X.Y., Xu, F.M., Zhang H.: L1/2 Regularization: A thresholding representation theory and a fast solver, \textit{IEEE Transactions on Neural Networks and Learning Systems},
2012, 24(7), pp.~1013--1027


\bibitem{Cao12}
Cao, W.F., Sun, J., Xu, Z.B.: Fast image deconvolution using closed-form thresholding formulas of $L_{q}(q=\frac{1}{2},\frac{2}{3})$ regularization,
\textit{Journal of Visual Communication and Image Representation}, 2013, 24(1), pp.~31--41


\bibitem{Zuo13}
Zuo, W.M., Meng, D.Y., Zhang, L., Feng, X.C., Zhang, D.: A generalized iterated shrinkage algorithm for non-convex sparse coding,
\textit{2013 IEEE International Conference on Computer Vision}, 2013, pp.~217--224


\bibitem{Cui14}
Cui, A.G., Peng, J.G., Li, H.Y., Zhang, C.Y., Yu, Y.C.: Affine matrix rank minimization problem via non-convex fraction function penalty, \textit{Journal of Computational and Applied Mathematics}, 2018, 336, pp.~353--374

\bibitem{Cui14jia}
Cui, A.G., Peng, J.G., Li, H.Y., Zhang, C.Y., Yu, Y.C.: Corrigendum to ``Affine matrix rank minimization problem via non-convex fraction function penalty" [J. Comput. Appl. Math. 336 (2018) 353-374], \textit{Journal of Computational and Applied Mathematics}, 2019, 352, pp.~478--485



\bibitem{Cui15}
Cui, A.G., Peng, J.G., Li H.Y., Wem M.: Nonconvex fraction function recovery sparse signal by convex optimization algorithm, 
\textit{arXiv:1905.05436v2}, 2019


\bibitem{don16}
Goldfarb, D., Ma, S.Q.: Convergence of fixed-point continuation algorithm for matrix rank minimization, \textit{Foundations of Computational Mathematics}, 2011, 11, pp.~183--210

\end{thebibliography}


\end{document}